\documentclass[12pt]{article}
\newtheorem{theorem}{Theorem}
\newtheorem{lemma}{Lemma}

\begin{document}

\title{Survival Probabilities for $N$-ary Subtrees on a
Galton-Watson Family Tree}

\author{\bf Ljuben R. Mutafchiev \footnote{Partial support is given
by the National Science Fund of the Bulgarian Ministry of
Education and Science, Grant No. VU-MI-105/2005.}
\\ American University in Bulgaria
\\ 2700 Blagoevgrad, Bulgaria and
\\ Institute of Mathematics and Informatics of the Bulgarian Academy
of Sciences
\\ \tt {ljuben@aubg.bg}}

\date{ }

\renewcommand{\baselinestretch}{1.3}

\maketitle

\begin{abstract}

The family tree of a Galton-Watson branching process may contain
$N$-ary subtrees, i.e. subtrees whose vertices have at least $N\ge
1$ children. For family trees without infinite $N$-ary subtrees,
we study how fast $N$-ary subtrees of height $t$ disappear as
$t\to\infty$.

\end{abstract}

\hspace {0.5cm} {\it Keywords:} Branching process; Family tree;
$N$-ary tree; Binary tree; Survival probability

\section{Introduction, Statement of the Results and Related Studies}

The family tree associated with a Bieneim\'{e}-Galton-Watson
process describes the evolution of a population in which each
individual, independently of the others, creates $k$ new
individuals with probability $p_k$ ($k=0,1,...$). We assume that
at generation zero there is single ancestor, called root of the
tree and let
\begin{equation}
f(s)=\sum_{k=0}^\infty p_k s^k
\end{equation}
to denote the probability generating function (pgf) of the
offspring distribution (with the convention that $f(1)=1$). We
recall the well known construction of a Galton-Watson family tree
noting that the individuals participating in the process and the
parent-child relations between them define the vertex-set and
arc-set of the tree, respectively (for more details see e.g.
Harris (1963, Ch. 6). By $\{Z_t, t=0,1,...\}$ we denote the
generation size process, defined by the following recurrence:
\begin{equation}
Z_0=1,  Z_{t+1}=\left\{\begin{array}{ll} X_1+...+X_{Z_t} &
\qquad\mbox{if}\qquad Z_t>0, \\ 0 & \qquad\mbox{if}\qquad Z_t=0,
\end{array} \right.
t=1,2,....
\end{equation}
Here $\{X_j, j=1,2,...\}$ are independent copies of a random
variable whose pgf is given by (1). In terms of trees $\{Z_t,
t=0,1,...\}$ present the sizes of the strata on the Galton-Watson
family tree. That is, $Z_t$ equals the number of vertices which
are at distance $t$ from the root. (We recall that the distance
between two vertices in a tree is determined by the number of arcs
in the path between them.)

The probabilities $P(Z_t>0)$ are often called survival
probabilities of the process. Their asymptotic behavior, as
$t\to\infty$, has been studied long ago by several authors. Let
\begin{equation}
\gamma_1=\lim_{t\to\infty}P(Z_t=0).
\end{equation}
It turns out that the parameters
\begin{equation}
a_1=f^\prime(\gamma_1), b_1=f^{\prime\prime}(\gamma_1)
\end{equation}
play important role in the study of survival probabilities of the
process. (If $\gamma_1=1$, here and further on, by $f^\prime(1)$,
$f^{\prime\prime}(1)$ and $f^{\prime\prime\prime}(1)$ we denote
the left derivatives of the power series (1) at the point $1$;
also note that $f^\prime(1)$ is the mean value of the offspring
distribution and $f^{\prime\prime}(1)$ is its second factorial
moment.) The following two results are well known and valid for
offspring distributions satisfying the inequality $p_0+p_1<1$.

{\it Result 1. [See Harris (1963, Ch. 1, Thms. 6.1.and 8.4).]} (i)
If $f^\prime(1)>1$ and $\gamma_1>0$, then $0<a_1<1$. (ii) For
$0<a_1<1$,
\begin{equation}
P(Z_t>0)=1-\gamma_1+d_1a_1^t+O(a_1^{2t}),
\end{equation}
as $t\to\infty$, where $d_1>0$ is certain constant.

{\it Result 2. [See Harris (1963, Ch. 1, Thm. 6.1 and Sect.
10.2).]} (i) If $f^\prime(1)=1$, then $\gamma_1=1$. (ii) If
$f^\prime(1)=1$ and $f^{\prime\prime\prime}(1)<\infty$, then
\begin{equation}
P(Z_t>0)\sim\frac{2}{b_1t}, t\to\infty.
\end{equation}

{\it Remark 1.} Result 2 was first obtained by Kolmogorov (1938).
It is also valid if $f^{\prime\prime}(1)<\infty$; see e.g.
Sevast'yanov (1971, Ch. 2, Sect. 2).

We will study special kinds of subtrees of a Galton-Watson family
tree. We will consider only rooted subtrees and call two such
subtrees disjoint if they do not have a common vertex different
from the root. Next, for fixed integer $N\ge 1$, we define a
complete infinite $N$-ary tree to be the family tree of a
deterministic branching process with offspring pgf $f(s)=s^N$. For
a branching process $\{Z_t, t=0,1,...\}$ defined by (2), we
introduce the random variable $V_N$, equal to the number of
complete disjoint and infinite $N$-ary subtrees rooted at the
ancestor. If we restrict the process up to its $t$th generation
($t=0,1,...$), or, equivalently, if we assume that the
Galton-Watson family tree is cut off at height t, we can similarly
define the random variable $V_{N,t}$ to be the count of the
complete disjoint $N$-ary subtrees of height at least $t$, which
are rooted at the ancestor. It is clear that
$V_N=\lim_{t\to\infty}V_{N,t}$ with probability $1$.

Further on we will assume that the offspring distribution $\{p_k,
k=0,1,...\}$ is such that $p_k<1$ for all $k$ and $p_k>0$ for some
$k>N$. For $N\ge 1$, we also let
\begin{equation}
\gamma_{N,0}=0, \gamma_{N,t}=P(V_{N,t}=0), t=1,2,...,
\end{equation}
and
\begin{equation}
\gamma_N=P(V_N=0)=\lim_{t\to\infty}\gamma_{N,t}.
\end{equation}
The last limit exists since the probabilities $\gamma_{N,t}$
monotonically increase in $t$.

In particular, for $N=1$, the event $\{V_1>0\}$ implies that the
family tree contains an infinite unary subtree (infinite path),
which means that the generations of the Galton-Watson process
never die. In the same way, event $\{V_2>0\}$ can be interpreted
as the set of trajectories of the process whose family trees grow
faster than binary splitting.

Another important observation follows from the well known
extinction criterion; Harris (1963, Ch. 1, Thm. 6.1). We give it
here in terms of the pgf (1) and probabilities (8) as follows: a
necessary and sufficient condition for $\gamma_1=1$ is
$f^\prime(1)\le 1$; if $f^\prime(1)>1$, then $\gamma_1$ is the
unique solution in $[0,1]$ of the equation
\begin{equation}
s=f(s).
\end{equation}
This enables one to restate Results 1(ii) and 2(ii) in terms of
counts of unary subtrees (recall (3) - (6) and Remark 1).

{\it Result $1^\prime$.} If $a_1<1$, then
 $$
 P(V_{1,t}>0\mid V_1=0)=d_1 a_1^t+O(a_1^{2t})
 $$
 as $t\to\infty$, where $d_1>0$ is certain constant.

{\it Result $2^\prime$.} If $a_1=1$ and
$f^{\prime\prime}(1)<\infty$, then
 $$
P(V_{1,t}>0\mid V_1=0)\sim\frac{2}{b_1 t}, t\to\infty.
 $$

 The main purpose of this present note is to study the survival of
 complete $N$-ary subtrees on a Galton-Watson family tree. We
 extend Results $1^\prime$ and $2^\prime$ to integer values of
 $N$ greater than $1$. Below we give the brief history of this
 problem.

 The question how to compute the probability that the
 Galton-Watson process possesses "the binary splitting property"
 was first raised, settled and solved by Dekking (1991). The
 general ($N\ge 2$) case was subsequently investigated by Pakes
 and Dekking (1991), who showed that the probability $\gamma_N$,
 defined by (8), is the smallest solution in $[0,1]$ of the
 equation
 \begin{equation}
 s=g_N(s),
 \end{equation}
 where
 \begin{equation}
 g_N(s)=\sum_{j=0}^{N-1}(1-s)^j f^{(j)}(s)/j!
 \end{equation}
 and $f(s)$ is the offspring pgf (1). We also point out that a
 particular case arising from a study of Mandelbrot's percolation process
 was previously considered by Chayes et al. (1988). Their problem is
equivalent to finding a condition on $p$ for $\gamma_8<1$ when
$f(s)=(1-p+ps)^9$. Furthermore, note that $g_1(s)=f(s)$ and
 thus, for $N=1$, eq. (10) reduces to (9). For particular offspring
 distributions, Pakes and Dekking (1991) encountered the following
 phenomenon: if $N\ge 2$, then there is
 a critical value $m_N^c$ for the offspring mean $f^\prime(1)$
 such that $\gamma_N=1$ if $f^\prime(1)<m_N^c$ and $\gamma_N<1$ if
 $f^\prime(1)\ge m_N^c$. Let $\gamma_N^c$ be the critical
 probability obtained when $f^\prime(1)=m_N^c$. It  turns out,
 for instance, that if $N=2$ and the offspring distribution is
 geometric, then $m_2^c=4$ and $\gamma_2^c=.75$; for a Poisson
 offspring distribution the same parameters are: $m_2^c=3.3509,
 \gamma_2^c=.4648$. (Further numerical results in this direction
 can be found in Pakes and Dekking (1991) and Yanev and Mutafchiev
 (2006).) This phenomenon is qualitatively different from what
 happens for $N=1$ where the extinction probability $\gamma_1=1$
 if $f^\prime(1)\le m_1^c=1$, except for the trivial case where
 $f(s)=s$ and $\gamma_1<1$ if $f^\prime(1)>1$. The case $N\ge 2$
 seems to be studied surprisingly later than the classical one
 when $N=1$. In fact, first assertions
 of the fundamental theorem on the existence of infinite unary
 subtrees on a Galton-Watson family tree appeared about 120 - 150
 years earlier and the problem was definitely settled around 1930.
 For more historical details, see e.g. Harris (1963) and
 Sevastyanov (1971). Recently, Yanev and Mutafchiev (2006) derived
 the probability distributions of the random variables $V_{N,t}$
 and $V_N$. The result for $V_{N,t}$ is given in a form of
 recurrence. Furthermore, the
 expression for the probability distribution of $V_N$ turns out
 to be very simple: its
 probability mass function equals the difference between two
 particular neighbor partial sums of the Taylor's expansion of $f(1)$
 around the point $\gamma_N$.

 To state our main results in an appropriate form, we extend
 notations (4) to integer values of $N\ge 2$. We set
 \begin{equation}
 a_N=g_N^\prime(\gamma_N), b_N=g_N^{\prime\prime}(\gamma_N)
 \end{equation}
and also recall definitions (7), (8) and (11).

\begin{theorem} If $\gamma_N\in(0,1)$ is the smallest solution of
eq. (10), then,

(i) for $N\ge 2$, we have $a_N\le 1$.

(ii) If $a_N<1$, then
 $$
 P(V_{N,t}>0\mid V_N=0)=d_N a_N^t+O(a_N^{2t})
 $$
 as $t\to\infty$, where $d_N>0$ is certain constant.

 (iii) If $a_N=1$, then

 (iiia) $b_N>0$, and,

 (iiib) for $N\ge 2$ and finite $b_N$,
 $$
P(V_{N,t}>0\mid V_N=0)\sim\frac{2}{\gamma_N b_N t}, t\to\infty.
 $$
\end{theorem}

Our paper is organized as follows. The proofs of the results are
presented in next Section 2. We recall there some old and
classical methods used in the theory of branching processes.
Section 3 contains few numerical results for particular offspring
distributions.

We conclude our introduction with a remark on studies which are
closely related to our model.

{\it Remark 2.} Pakes and Dekking (1991) noticed that there are
links between complete infinite $N$-ary subtrees on a
Galton-Watson family tree, Mandelbrot's percolation process
studied by by Chayes et al. (1988) and results obtained by
Pemantle (1988) and related to a model of a reinforced random
walk. In particular, Pemantle (1988) established the following
criterion for $\gamma_N<1$ (see his Lemma 5 or Pakes and Dekking
(1991, pp. 356-357)): if for some $s_0\in(0,1)$ we have
$g_N(s_0)\le s_0$, then $\gamma_N\le s_0$. Here we also indicate a
relationship between the $N$-ary subtrees phenomenon and the
existence of a $k$-core in a random graph. The $k$-core of a graph
is the largest subgraph with minimum degree at least $k$. This
concept was introduced by Bollob\'{a}s (1984) in the context of
finding large $k$-connected subgraphs of random graphs. He
considered the Erd\"{o}s-R\'{e}nyi random graph $G(n,p)$ with $n$
vertices in which the possible arcs are present independently,
each with probability $p$. If we set $p=\lambda/n$, where
$\lambda>0$ is a constant, it is natural to ask: for $k\ge 3$,
what is the critical value $\lambda_c(k)$ of $\lambda$ above which
a (non-empty) $k$-core first appears in $G(n,\lambda/n)$ with
probability tending to $1$ as $n\to\infty$. To answer this
question Pittel et al. (1996) considered a Galton-Watson family
tree rooted at a vertex $x_0$ (ancestor) of the graph
$G(n,\lambda/n)$ and assume that the offspring distribution of the
branching process is Poisson with mean $\lambda$. Let $B_k$ denote
the event that $x_0$ has at least $k$ children each of which has
at least $k-1$ children each of which has at least $k-1$ children,
and so on. It is clear that this assumption slightly modifies the
concept of a complete infinite $(k-1)$-ary subtree (the only
difference occurs in the assumption for the offspring number of
the ancestor $x_0$). Pittel et al. (1996) found the threshold
$\lambda_c(k)$ for the emergence of a non-trivial $k$-core in
$G(n,\lambda/n)$ and showed that, except at the critical value,
the number of vertices in the $k$-core approaches $P(B_k)n$ as
$n\to\infty$. Their results also showed that a giant $k$-core
appears suddenly when the number of arcs in the random graph
reaches $c_k n/2$, where the constants $c_k$ are explicitly
computed. There is a remarkable coincidence between constants
$c_k$ and the critical means $m_{k-1}^c (k=3,4,5)$ of the Poisson
offspring distributions which yield existence of $(k-1)$-ary
subtree on a Galton-Watson family tree given by Yanev and
Mutafchiev (2006, p. 232). The idea of embedding a Poisson
branching process in the random graph model was recently developed
by Riordan (2007) who gave a new proof of the results of Pittel et
al. (1996) and extended them to a general model of inhomogeneous
random graphs with independence between their arcs.

\section{Proofs of the Results}

First, we recall Pakes and Dekking (1991) result: the probability
$\gamma_N$, defined by (8), is the smallest solution in $[0,1]$ of
eq. (10). To prove part (i) of the theorem, note that $\gamma_N>0$
implies that $g_N(0)>0$. Therefore, for $s\in[0,\gamma_N)$, the
graph of the function $y=g_N(s)$ lies above the diagonal of the
unit square in the coordinate system $sOy$. At $s=\gamma_N$ the
curve $y=g_N(s)$ crosses or touches the diagonal $y=s$. If it
touches it, then $a_N=g_N^\prime(\gamma_N)=1$. If $y=g_N(s)$
crosses the diagonal, then, for some sufficiently small
$\epsilon\in (0,\gamma_N)$, we have
$g_N(\gamma_N-\epsilon)\ge\gamma_N-\epsilon$ and
$g_N(\gamma_N+\epsilon)\le\gamma_N+\epsilon$. Hence
$g_N(\gamma_N+\epsilon)-g_N(\gamma_N-\epsilon)\le 2\epsilon$.
Therefore the derivative
 $$
 g_N^\prime(s)=(1-s)^{N-1}f^{(N)}(s)/(N-1)!
 $$
 should not exceed $1$ for certain
 $s=s_\epsilon\in(\gamma_N-\epsilon,\gamma_N+\epsilon)$, by the
 mean value theorem. Letting $\epsilon\to 0$ and using the
 continuity of $g_N^\prime(s)$, we get assertion (i).

 Assertion (ii) can be obtained using a result on iterations of
 functions due to Koenigs (1884) (see also Harris (1963, Ch. 1,
 Sect. 8.3)). Below we state a suitable modification of it as a
 separate lemma. The proof follows the same line of reasoning as
 in Harris (1963, Ch. 1, Thm. 8.4).

  \begin{lemma} Let
$$
 h(s)=\sum_{j=0}^\infty h_j s^j
 $$
 ($h_j$ real) be a function, which is analytic in $\mid s\mid<1$,
 strictly increasing in $[0,1]$ and such that $h(1)=1$. Let
 \begin{equation}
 h_0(s), h_1(s)=h(s), h_{t+1}=h(h_t(s)), t=1,2,...
 \end{equation}
 be the sequence of iterations of $h(s)$. Suppose that the
 equation
 \begin{equation}
 s=h(s)
 \end{equation}
 has a solution in $[0,1]$ and let $q$ be the least one in
 $[0,1]$. If $q$ satisfies $h^\prime(q)<1$, then
 $$
 h_t(0)=q-d[h^\prime(q)]^t+O([h^\prime(q)]^{2t})
 $$
 as $t\to\infty$, where $d>0$ denotes an absolute constant.
  \end{lemma}

  We will apply Lemma 1 setting $h(s)=g_N(s)$. Define the iterations
  $g_{N,t}(s)$ of the
  function $g_N(s)$ as in (13). Also, recall
  that $g_N(1)=1$ and $\gamma_{N,0}=0$ (see definitions (11) and
  (7), respectively). We
  set $q=\gamma_N$ in eq. (14). Then, we use the recurrence
  $\gamma_{N,t}=g_N(\gamma_{N,t-1})$;
  see Yanev and Mutafchiev (2006, p. 227). Iterating
  $t$ times as in (13), we get
  $\gamma_{N,t}=P(V_{N,t}=0)=g_{N,t}(0)$. Hence, by Lemma 1 and
  notation $(12_1)$,
  $$
  \gamma_{N,t}=\gamma_N-d_N^\prime a_N^t+O(a_N^{2t})
  $$
  as $t\to\infty$, where $d_N^\prime>0$ denotes an absolute constant.
  Dividing both sides of this equality by $\gamma_N$ and writing
  conditional probabilities for $V_{N,t}$, we obtain
  assertion (ii) with $d_N=d_N^\prime/\gamma_N$.

  To prove (iiia), let us assume that $b_N\le 0$ (see notation
  $(12_2)$). This shows that $g_N^\prime(s)$ decreases
  in a neighborhood of $s=\gamma_N$. Hence, there exists a
  sufficiently small number $\delta>0$ such that, for any
  $s\in(\gamma_N-\delta,\gamma_N]$, we have $g_N^\prime(s)\ge
  g_N^\prime(\gamma_N)=a_N=1$. Therefore, $[g_N(s)-s]^\prime\ge
  0$, and so, the function $g_N(s)-s$ increases in
  $(\gamma_N-\delta,\gamma_N]$. Thus, for any
  $s\in(\gamma_N-\delta,\gamma_N)$, we have $g_N(s)-s\le
  g_N(\gamma_N)-\gamma_N=0$. Combining the inequalities $g_N(s)\le
  s, g_N(0)>0$ and using the continuity of $g_N(s)$, we conclude
  that there is some $s_0<\gamma_N$ that solves eq. (10). This
  contradicts the assumption that $\gamma_N$ is the smallest
  solution in $(0,1]$ of eq. (10). So, (iiia) is proved.

  The asymptotic given in assertion (iiib) will also follow from
  classical results on iterations of analytic functions,
  increasing on a segment of the real axis. One possible proof may use
  a general result of Harris (1963, Ch. 1, Lemma 10.1)
  establishing uniform asymptotics for $1/[1-h_t(s)]$, where
  $h_t(s)$ denote the iterations defined by (13) and the complex variable
  $s$ varies in some particular subsets of the unit disc. In our
  case it suffices, however, to consider the behavior of $h_t(s)$
  only at $s=0$. The problem turns out to be similar to that for the
  critical branching process which was studied first by Kolmogorov
  (1938). The proof of the next lemma repeats the arguments given
  by Sevast'yanov (1971, Ch. 2, Sect. 2).

  \begin{lemma} Suppose that $h(s), h_t(s), t=0,1,...$ and $q$ are
  the same as in Lemma 1. Furthermore, suppose that
  $h^\prime(q)=1$ and $h^{\prime\prime}(q)\in (0,\infty)$. Then, we have
  \begin{equation}
  \frac{1}{q-h_t(0)}=\frac{th^{\prime\prime}(q)}{2}+o(t)
  \end{equation}
  as $t\to\infty$.
  \end{lemma}

  To show how assertion (iiib) follows from this lemma we set
  in both sides of (15): $q=\gamma_N, h(s)=g_N(s),
  h_t(0)=g_{N,t}(0)=\gamma_{N,t}$ and
  $h^{\prime\prime}(q)=g_N^{\prime\prime}(\gamma_N)=b_N$. Thus, we
  obtain
  \begin{equation}
  \frac{1}{\gamma_N-\gamma_{N,t}}=\frac{tb_N}{2}+o(t),
  t\to\infty.
  \end{equation}
  To complete the proof of (iiib) it remains to take the reciprocal
  of (16), divide both sides by $\gamma_N$ and convert the ratio
  $\gamma_{N,t}/\gamma_N$ into conditional probability for
  $V_{N,t}$.

 \section{Numerical Results}

 {\it Geometric distribution.} We look at the case, where
 $$
 f(s)=\frac{1-p}{1-ps},
 g_N(s)=1-\left[\frac{p(1-s)}{1-ps}\right]^N, 0<p<1.
 $$
Pakes and Dekking (1991) established in this case that the
critical mean for $N=2$ is $m_2^c=4$ which implies that the
critical value for the parameter $p$ is $p_2^c=4/5$. It is easy to
see that the least solution in $[0,1]$ of eq. (10) is
$\gamma_2^c=3/4$. Calculating the first two derivatives of
$g_2(s)$ at $s=3/4$, we get $a_2^c=1, b_2^c=2$. Therefore, by
assertion (iiib) of Theorem 1,
 $$
 P(V_{2,t}>0\mid V_2=0)\sim\frac{4}{3t}, t\to\infty.
 $$

 {\it Poisson distribution.} The Poisson offspring distribution
 has the pgf $f(s)=e^{m(s-1)}, m>0$. Whence
 $$
 g_N(s)=e^{m(s-1)}\sum_{j=0}^{N-1}[(1-s)m]^j/j!.
 $$
 In this case $m_2^c=3.3509$ and the least solution in $[0,1]$ of
 eq. (10) is $\gamma_2^c=.4648$ (see Yanev and Mutafchiev (2006)).
 Numerical computations with greater level of accuracy show
 that $a_2^c=1, b_2^c=1.48235$ and by Theorem 1(iiib),
 $$
 P(V_{2,t}>0\mid V_2=0)\sim\frac{2.9028}{t},
 t\to\infty.
 $$

 {\it One-or-many distribution.} This is a two-parameter family of
 discrete distributions defined for some $p\in(0,1)$ and integer
 $r>N>1$ by the equalities: $p_r=p, p_1=1-p$. Clearly,
 $f(s)=(1-p)s+ps^r$, and hence
 $$
 g_N(s)=1-p\sum_{j=N}^r {r\choose j}(1-s)^j s^{r-j}.
 $$
 Pakes and Dekking (1991) showed that if $r=N+1$, then
 $\gamma_N^c=1/N^2$ and the threshold value of the parameter $p$
 is $p_N^c=(1-1/N)(1-1/N^2)^{-N}$. For $r=3$ and $N=2$, we have
 $g_2^\prime(s)=6ps(1-s), p_2^c=8/9, \gamma_2^c=1/4$. Thus, we get
 $a_2^c=1, b_2^c=8/3$, and hence by Theorem 1(iiib),
 $$
P(V_{2,t}>0\mid V_2=0)=\frac{3}{t}, t\to\infty.
 $$

\end{document}